\catcode`\^^Z=9
\catcode`\^^M=10
\output={\if N\header\headline={\hfill}\fi
\plainoutput\global\let\header=Y}
\magnification\magstep1
\tolerance = 500
\hsize=14.4true cm
\vsize=22.5true cm
\parindent=6true mm\overfullrule=2pt
\newcount\kapnum \kapnum=0
\newcount\parnum \parnum=0
\newcount\procnum \procnum=0
\newcount\nicknum \nicknum=1
\font\ninett=cmtt9

\font\ninebf=cmbx9

\font\sixbf=cmbx6
\font\ninesl=cmsl9

\font\nineit=cmti9

\font\ninerm=cmr9

\font\sixrm=cmr6
\font\ninei=cmmi9
\font\eighti=cmmi8
\font\sixi=cmmi6
\skewchar\ninei='177 \skewchar\eighti='177 \skewchar\sixi='177
\font\ninesy=cmsy9
\font\eightsy=cmsy8
\font\sixsy=cmsy6
\skewchar\ninesy='60 \skewchar\eightsy='60 \skewchar\sixsy='60
\font\titelfont=cmr10 scaled 1440
\font\paragratit=cmbx10 scaled 1200

\font\name=cmcsc10
\font\emph=cmbxti10

\font\tenmsbm=msbm10
\font\sevenmsbm=msbm7
%

%

%
\font\teneufm=eufm10
\font\seveneufm=eufm7
\font\fiveeufm=eufm5
\newfam\eufmfam
\textfont\eufmfam=\teneufm
\scriptfont\eufmfam=\seveneufm
\scriptscriptfont\eufmfam=\fiveeufm

\font\tenmsam=msam10
\font\sevenmsam=msam7
\font\fivemsam=msam5
\newfam\msamfam
\textfont\msamfam=\tenmsam
\scriptfont\msamfam=\sevenmsam
\scriptscriptfont\msamfam=\fivemsam
\font\tenmsbm=msbm10
\font\sevenmsbm=msbm7
\font\fivemsbm=msbm5
\newfam\msbmfam
\textfont\msbmfam=\tenmsbm
\scriptfont\msbmfam=\sevenmsbm
\scriptscriptfont\msbmfam=\fivemsbm
\def\Bbb#1{{\fam\msbmfam\relax#1}}
\def\cz{{\kern0.4pt\Bbb C\kern0.7pt}
}
\def\ez{{\kern0.4pt\Bbb E\kern0.7pt}
}
\def\fz{{\kern0.4pt\Bbb F\kern0.3pt}}
\def\gz{{\kern0.4pt\Bbb Z\kern0.7pt}}
\def\hz{{\kern0.4pt\Bbb H\kern0.7pt}
}
\def\kz{{\kern0.4pt\Bbb K\kern0.7pt}
}
\def\nz{{\kern0.4pt\Bbb N\kern0.7pt}
}
\def\oz{{\kern0.4pt\Bbb O\kern0.7pt}
}
\def\rz{{\kern0.4pt\Bbb R\kern0.7pt}
}
\def\sz{{\kern0.4pt\Bbb S\kern0.7pt}
}
\def\pz{{\kern0.4pt\Bbb P\kern0.7pt}
}
\def\qz{{\kern0.4pt\Bbb Q\kern0.7pt}
}
\newskip\ttglue
\def\ninepoint{\def\rm{\fam0\ninerm}%
  \textfont0=\ninerm \scriptfont0=\sixrm \scriptscriptfont0=\fiverm
  \textfont1=\ninei \scriptfont1=\sixi \scriptscriptfont1=\fivei
  \textfont2=\ninesy \scriptfont2=\sixsy \scriptscriptfont2=\fivesy
  \textfont3=\tenex \scriptfont3=\tenex \scriptscriptfont3=\tenex
  \def\it{\fam\itfam\nineit}%
  \textfont\itfam=\nineit
  \def\sl{\fam\slfam\ninesl}%
  \textfont\slfam=\ninesl
  \def\bf{\fam\bffam\ninebf}%
  \textfont\bffam=\ninebf \scriptfont\bffam=\sixbf
   \scriptscriptfont\bffam=\fivebf
  \def\tt{\fam\ttfam\ninett}%
  \textfont\ttfam=\ninett
  \tt \ttglue=.5em plus.25em minus.15em
  \normalbaselineskip=11pt
  \font\name=cmcsc9
  \let\sc=\sevenrm
  \let\big=\ninebig
  \setbox\strutbox=\hbox{\vrule height8pt depth3pt width0pt}%
  \normalbaselines\rm
  \def\sl{\it}}

\headline={\ifodd\pageno\rightheadline\else\leftheadline\fi}
\def\rightheadline{\ninepoint Paragraphen"uberschrift\hfill\folio}
\def\leftheadline{\ninepoint\folio\hfill Chapter"uberschrift}
\let\header=Y
\def\titel#1{\need 9cm \vskip 2truecm
\parnum=0\global\advance \kapnum by 1
{\baselineskip=16pt\lineskip=16pt\rightskip0pt
plus4em\spaceskip.3333em\xspaceskip.5em\pretolerance=10000\noindent
\titelfont Chapter \uppercase\expandafter{\romannumeral\kapnum}.
#1\vskip2true cm}\def\leftheadline{\ninepoint
\folio\hfill Chapter \uppercase\expandafter{\romannumeral\kapnum}.
#1}\let\header=N
}
\def\Titel#1{\need 9cm \vskip 2truecm
\global\advance \kapnum by 1
{\baselineskip=16pt\lineskip=16pt\rightskip0pt
plus4em\spaceskip.3333em\xspaceskip.5em\pretolerance=10000\noindent
\titelfont\uppercase\expandafter{\romannumeral\kapnum}.
#1\vskip2true cm}\def\leftheadline{\ninepoint
\folio\hfill\uppercase\expandafter{\romannumeral\kapnum}.
#1}\let\header=N
}
\def\need#1cm {\par\dimen0=\pagetotal\ifdim\dimen0<\vsize
\global\advance\dimen0by#1 true cm
\ifdim\dimen0>\vsize\vfil\eject\noindent\fi\fi}
\def\neupara#1{\par\penalty-2000
\procnum=0\global\advance\parnum by 1
\vskip1cm\noindent{\paragratit \the\parnum. #1}%
\def\rightheadline{\ninepoint\S\the\parnum.\ #1\hfill \folio}%
\vskip 8mm\noindent}
\def\Proclaim #1 #2\finishproclaim {\bigbreak\noindent
{\bf#1\unskip{}. }{\it#2}\medbreak\noindent}
%
\gdef\proclaim #1 #2 #3\finishproclaim {\bigbreak\noindent%
\global\advance\procnum by 1
{%
{\relax\ifodd \nicknum
\hbox to 0pt{\vrule depth 0pt height0pt width\hsize
   \quad \ninett#3\hss}\else {}\fi}%
\bf\the\parnum.\the\procnum\ #1\unskip{}. }
{\it#2}
\immediate\write\num{\string\def
 \expandafter\string\csname#3\endcsname
 {\the\parnum.\the\procnum}}
\medbreak\noindent}
\newcount\stunde \newcount\minute \newcount\hilfsvar
\def\uhrzeit{
    \stunde=\the\time \divide \stunde by 60
    \minute=\the\time
    \hilfsvar=\stunde \multiply \hilfsvar by 60
    \advance \minute by -\hilfsvar
    \ifnum\the\stunde<10
    \ifnum\the\minute<10
    0\the\stunde:0\the\minute~Uhr
    \else
    0\the\stunde:\the\minute~Uhr
    \fi
    \else
    \ifnum\the\minute<10
    \the\stunde:0\the\minute~Uhr
    \else
    \the\stunde:\the\minute~Uhr
    \fi
    \fi
    }
 \def\calB{{\cal B}}

\def\calE{{\cal E}}

\def\calM{{\cal M}} \def\calN{{\cal N}}
 
 \def\calR{{\cal R}}
\def\calS{{\cal S}} 
\def\calU{{\cal U}} 
 
 \def\calZ{{\cal Z}}

\def\End{\mathop{\rm End}\nolimits}

\def\GL{\mathop{\rm GL}\nolimits}

\def\id{\mathop{\rm id}\nolimits}

\def\kernel{\mathop{\rm kernel}\nolimits}

\def\U{{\rm U}}

\def\boxit#1{
  \vbox{\hrule\hbox{\vrule\kern6pt
  \vbox{\kern8pt#1\kern8pt}\kern6pt\vrule}\hrule}}
\def\Boxit#1{
  \vbox{\hrule\hbox{\vrule\kern2pt
  \vbox{\kern2pt#1\kern2pt}\kern2pt\vrule}\hrule}}

\def\smallni{\smallskip\noindent }
\def\medni{\medskip\noindent }

\def\Isom{\mathop{\;{\buildrel \sim\over\longrightarrow }\;}}
\def\lo{\longrightarrow}

\def\loma{\longmapsto}
\def\betr#1{\vert#1\vert}
\def\spitz#1{\langle#1\rangle}

\def\pii{\pi {\rm i}}

\def\square{\hbox{\hbox to 0pt{$\sqcup$\hss}\hbox{$\sqcap$}}}
\def\qed{\ifmmode\square\else{\unskip\nobreak\hfil
\penalty50\hskip3em\null\nobreak\hfil\square
\parfillskip=0pt\finalhyphendemerits=0\endgraf}\fi}
\def\pn{\the\parnum.\the\procnum}
\def\downmapsto{{\buildrel
        {\vbox{\hbox{\hskip.2pt$\scriptstyle-$}}}
        \over{\raise7pt\vbox{\vskip-4pt\hbox{$\textstyle\downarrow$}}}}}

\def\NisoN{1.4}

\def\DefAu{2.2}
\def\sigK{2.3}

\def\CalM{3.1}
\def\RamLoc{3.2}
\def\GiSi{4.1}
\def\DeRel{5.1}
\def\KlaDet{5.2}
\def\ProInt{5.3}

\def\MnotN{5.6}
\def\MT{5.7}
\nopagenumbers
\immediate\newwrite\num
\nicknum=0  
\let\header=N
\def\transpose#1{\kern1pt{^t\kern-1pt#1}}%

\immediate\openout\num=calabi8.num
\immediate\newwrite\num\immediate\openout\num=ball-vector.num
\def\RAND#1{\hbox to 0mm{\hss\vtop to 0pt{%
  \raggedright\ninepoint\parindent=0pt%
  \baselineskip=1pt\hsize=2cm #1\vss}}\noindent}
\noindent
\centerline{\titelfont Vector valued  modular forms on  three dimensional ball}%
\def\leftheadline{\ninepoint\folio\hfill
 Vector valued  modular forms on  three dimensional ball}%
\def\rightheadline{\ninepoint Introduction\hfill \folio}%
\headline={\ifodd\pageno\rightheadline\else\leftheadline\fi}
\vskip 1.5cm
\leftline{\it \hbox to 6cm{Eberhard Freitag\hss}
Riccardo Salvati  
Manni  }
  \leftline {\it  \hbox to 6cm{Mathematisches Institut\hss}
Dipartimento di Matematica, }
\leftline {\it  \hbox to 6cm{Im Neuenheimer Feld 288\hss}
Piazzale Aldo Moro, 2}
\leftline {\it  \hbox to 6cm{D69120 Heidelberg\hss}
 I--00185 Roma, Italy. }
\leftline {\tt \hbox to 6cm{freitag@mathi.uni-heidelberg.de\hss}
salvati@mat.uniroma1.it}
\vskip1cm
\centerline{\paragratit \rm  2014}%
\vskip5mm\noindent%
\let\header=N%
{\paragratit Introduction}%
\medni
In the paper [CG],  Cl\'ery and van der Geer determined generators for some modules of vector
valued Picard modular forms on the two dimensional ball. In this paper we consider the case
of a three dimensional ball with the action of the Picard modular group $\Gamma_3[\sqrt{-3}]$
(see Sect.~3).  The corresponding
modular variety of dimension 3 is a copy of the Segre cubic.
\smallskip
Vector valued Picard modular forms on the $n$-ball $\calB_n$
belong to rational representations of the complexification
or the maximal compact group of the unitary group $\U(1,n)$, which is the group
$\GL(1,\cz)\times\GL(n,\cz)$. Here we consider the representation
$$\varrho_r(k_1,k_2)=k_1^{r}k_2\qquad(r\in\gz).$$
A similar representation in a Siegel case has been treated in [FS2].
We denote by $\calM(r)$ the space of modular forms $f:\calB_n\to\cz^n$ which belong
to this representation. The direct sum
$$\calM=\bigoplus_{r\in\gz}\calM(r)$$
is a module over the ring of scalar-valued modular forms. 
\smallskip
In the case $\Gamma_3[\sqrt{-3}]$
this ring is generated by 5 forms $T_1,\dots,T_5$ of weight 3 which satisfy the relation
of a Segre cubic, [FS1,Ko].
We will determine the structure
of the module $\calM$. For this we consider the submodule $\calN$ of $\calM$, generated
by 10 Cohen-Rankin brackets $\{T_i,T_j\}$. They are elements of $\calM(5)$.
One of our main results is that $\calM$ and $\calN$ nearly agree. They differ only in the 
lowest possible degrees $r=5$ and $r=8$. (We always have $r\equiv 2$ mod $3$ if
$\calM(r)$ is not zero.)
An extra form in weight $5$ will be constructed explicitely. This form and those in    $\calN$   generate 
$\calM$.
\smallskip
To get a proof,  we first   determine the structure of $\calN$. There 
are some obvious relations
between the Cohen-Rankin brackets and  also the Segre relation induces a relation between
them. That these simple relations are defining relations 
(see Proposition \DeRel) rests on
a pure algebraic statement about differential modules which is developed in Sect.~1.
In Sect.~2 we develop the   framework for vector valued ball modular forms and in Sect.~3
we describe the group of our interest $\Gamma_3[\sqrt{-3}]$, 
the congruence group of level $\sqrt{-3}$
in the full Picard modular group with respect to $\qz(\sqrt{-3})$. We describe its ring of
modular forms, the relation to the Segre cubic and the ramification locus.
\smallskip
In Sect.~4 we study some
special modular forms which are related to the tangent bundle of the Segre cubic.
They are needed for the prove of the basic relation between $\calM$ and its submodule
$\calN$,
$$\calM=\bigcap {1\over T_i^2}\calN$$
which will given in  Sect.~5. The  structure theorem for $\calM$ can be derived  from this result.
Some computer calculations are necessary.
\smallskip
In our main result, Theorem \MT, we give generators of the module $\calM$ and
we produce the Hilbert functions of the modules $\calN$, $\calM$. 
\neupara{Differential modules over graded algebras}%
Let
$A=\bigoplus_{d=0}^\infty A_d$ be a finitely generated graded algebra over
a field $K=A_0$ of characteristic 0. We assume that $A$ is an integral
domain and denote its field of fractions by $Q(A)$.
We consider the K\"ahler differential module
$$\Omega=\Omega(Q(A)/K).$$
Recall that this is a $Q(A)$-vector space together with a
$K$-linear derivation $d:Q(A)\to \Omega$. 
The dimension of $\Omega$ equals the transcendental degree of $Q(A)$ and $\Omega$
is generated by the image of $d$.
In the following, we denote by $\deg(f)$ the degree of
a non-zero homogeneous element of $A$. 
For two non-zero homogeneous elements of positive degree $f,g\in A$ we define
$$\{f,g\}:=\deg(g)g df-\deg(f)fdg.$$
Another way to write this is
$$\{f,g\}={g^{\deg(f)+1}\over f^{\deg(g)-1}}\; d\Bigl(
{f^{\deg(g)}\over g^{\deg(f)}}\Bigr).$$
This is a skew-symmetric $K$-bilinear pairing and it satisfies
the following rule
$$\deg(h)h\{f,g\}=\deg(g)g\{f,h\}+\deg(f)f\{h,g\}.$$
\proclaim
{Definition}
{We denote by $\calN$ the $A$-module that is generated
by all $\{f,g\}$ where $f,g$ are homogeneous elements of positive degree in $A$.}
defN%
\finishproclaim
We are interested in a finite presentation of $\calN$.
There is no difficulty to get a finite system of generators.
Let $A=K[f_1,\dots, f_m]$, ($f_i$ homogenous). Then
$\{f_i,f_j\}$ are generators of $\calN$. It is more involved to get
defining relations.
\smallskip
We use the notation $d_i=\deg(f_i)$. A polynomial
$P\in K[X_1,\dots,X_m]$ is called isobaric of weight $k$
(with respect to
$(d_1,\dots,d_m)$) if it is of the form
$$P=\sum_{d_1\nu_1+\cdots+d_m\nu_m=k}a_{\nu_1,\dots,\nu_m}
X_1^{\nu_1}\cdots X_m^{\nu_m}.$$
Then the Euler relation
$$\sum_{\nu=1}^m d_\nu{\partial P\over\partial X_\nu} X_\nu=kP$$
holds.
\smallskip
The ideal of relations between $f_1,\dots,f_m$ is
generated by isobaric polynomials. Let $R(f_1,\dots,f_m)=0$ be
an isobaric relation.  Differentiation gives
$$\sum_{\nu=1}^m(\partial_\nu R)df_\nu=0\quad\hbox{where}\quad
\partial_\nu R:={\partial R\over\partial X_\nu}(f_1,\dots,f_m).$$
From this relation and the Euler relation we derive
$$\sum_{\nu=1}^m(\partial_\nu R)\{f_\nu,f_\mu\}=0\quad (\mu\ \hbox{arbitrary}).$$
We want to formalize this and introduce a module $\calN'$ which
is defined by the so far known relations.
\proclaim
{Definition}
{
We denote by $\calN'$ the $A$-module that is generated by symbols
$[f_i,f_j]$ with the following defining relations:
$$d_kf_k[f_i,f_j]=d_jf_j[f_i,f_k]+d_if_i[f_k,f_j],\quad [f_i,f_j]+[f_j,f_i]=0.
\leqno\hbox{\rm(1)}$$
For each isobaric relation $R$  between the $f_1,\dots,f_m$ one has
$$\sum_{\nu=1}^m(\partial_\nu R)[f_\nu,f_\mu]=0\quad
(\mu\ \hbox{arbitrary}).
\leqno\hbox{\rm(2)}$$
}
Nstrich%
\finishproclaim
It is of course enough to take for $R$ a system of generators of the ideal
of all relations.

\smallskip
There is a natural surjective homomorphism
$$\calN'\lo\calN,\quad [f_i,f_j]\loma \{f_i,f_j\}.$$
We notice that $\calN$ is torsion free for trivial reasons, but it is not clear that
$\calN'$ is   torsion free too.
\smallskip
Under certain circumstances, $\calN'\to\calN$ is an isomorphism.
To work this out, we consider an arbitrary relation in $\calN$
$$\sum_{i<j}P_{ij}\{f_i,f_j\}=0,\quad P_{ij}\in A.$$
We multiply this relation by $d_1f_1$ and insert
$$d_1f_1\{f_i,f_j\}=d_if_i\{f_1,f_j\}-d_jf_j\{f_1,f_i\}.$$
Then we obtain the relation
$$\sum_{j} P_{j}\{f_1,f_j\}=0,$$
where the elements $P_j\in A$ are defined as
$$P_j=\sum_{i<j} d_if_i P_{ij}-\sum_{i>j}d_if_i P_{ji}.$$
Let $n$ be the transcendental degree of $Q(A)$. We can assume that
$f_1,\dots,f_n$ are independent. Then each $f_k$, $k>n$, satisfies
an algebraic relation
$$R_k(f_1,\dots,f_n,f_k)=0.$$
Here $R_k$ is an irreducible polynomial in the variables
$X_1,\dots,X_n,X_k$.
Now we make use of the relation
$$(\partial_k R_k)\{f_1,f_k\}+\sum_{\nu=1}^n(\partial_\nu R_k)\{f_1,f_\nu\}=0.$$
We have to use the elements (from the ring $A$)
$$\Pi:=\prod_{k=n+1}^m \partial_k R_k,\quad \Pi^{(k)}:={\Pi\over \partial_k P_k}.$$
We multiply the original relation  by $\Pi$:
$$\Pi\sum_{j} P_{j}\{f_1,f_j\}=0.$$
For $k>n$ we have the formula
$$\Pi\{f_1,f_k\}=\Pi^{(k)}(\partial_kR_k)\{f_1,f_k\}=-
\Pi^{(k)}\sum_{j=1}^n(\partial_j R_k)\{f_1,f_j\}.$$
Now we can eliminate the $\{f_1,f_k\}$ for $k>n$ to produce
a relation between the $\{f_1,f_i\}$, $2\le i\le n$. But these elements are
independent. Hence the coefficients of the relation must vanish.
A simple calculation now gives the following lemma.
\proclaim
{Lemma}
{Let
$$\sum_{i<j}P_{ij}\{f_i,f_j\}=0,\quad P_{ij}\in A.$$
Then the elements
$$P_j=\sum_{i<j} d_if_i P_{ij}-\sum_{i>j}d_if_i P_{ji}$$
satisfy the following system of relations.
$$P_j\Pi=\sum_{k=n+1}^m (\partial_j R_k)P_k\Pi^{(k)}\quad (1\le j\le n).$$
\smallni
{\bf Supplement.}
Conversely,  these relations imply in $\calN'$ the relation
$$f_1\Pi\sum_{i<j}P_{ij}[f_i,f_j]=0.$$
}
SysR%
\finishproclaim
For the proof of the supplement we just have to notice that the calculations above
only use the defining relations of $\calN'$.\qed
\smallskip
Let us assume that multiplication by $f_1\Pi$ is injective on $\calN'$.
Then we see that $\sum P_{ij}\{f_i,f_j\}=0$ implies 
$\sum P_{ij}[f_i,f_j]=0$. Hence $\calN'\to \calN$ is an isomorphism
and $\calN'$ must be torsion free.
This gives the following result.
\proclaim
{Proposition}
{
Assume that the
$f_1,\dots,f_n$ is a transcendental basis such that  each $f_k$, 
$n<k\le m$, satisfies
an irreducible algebraic relation
$$R_k(f_1,\dots,f_n,f_k)=0.$$
The homomorphism $\calN'\to\calN$ is an isomorphism
if and only if $\calN'$ is torsion free. For this it suffices that
multiplications by $f_1$ and $\partial_k R_k$ ($n<k\le m$) are injective
on $\calN'$.}
NisoN%
\finishproclaim
\neupara{The extended ball}%
Let $V$ be a complex vector space of dimension $n+1$ and let $\spitz{\cdot,\cdot}$
be a non degenerated hermitian form
of signature $(1,n)$.
We consider the projective space $\pz(V)=(V-\{0\})/\cz^*$ and the natural projection
$$V-\{0\}\lo \pz(V),\quad v\loma [v].$$
Let
$$\tilde\calB:=\{v\in V;\quad \spitz{v,v}>0\}$$
be the set of all vectors of positive norm $\spitz{v,v}>0$ and  $\calB$ its image
in the projective space. This is a model of the complex $n$-ball.
The unitary group $\U(V)$ acts on $\calB$ and on $\tilde\calB$.
\smallskip
We choose a vector
$e\in V$ with positive norm $\spitz{e,e}>0$ and we consider
the orthogonal complement
$\calZ=e^\perp$ which is a negative definite space of dimension $n$.
We have $V=\cz e\oplus \calZ$. Sometimes we write the elements $v\in V$ in the form
$$v=Ce+z=\pmatrix{C\cr z}.$$
Then we can write the elements of $\End(V)$ as matrices
$$p=\pmatrix{a&b\cr c&d},\quad a\in\cz,\ b\in \calZ^*,\ c\in \calZ,\ d\in\End(\calZ),$$
such that the action on $V=\cz e+\calZ$ is given by
$$\pmatrix{a&b\cr c&d}\pmatrix{C\cr z}=\pmatrix{aC+b(z)\cr Cc+d(z)}.$$
For the multiplication of two of such matrices one has to make use of the
canonical isomorphism $\calZ\otimes \calZ^*\to\End(\calZ)$.
\smallskip
We denote by
$$\calB_\calZ:=\{z\in \calZ;\quad -\spitz{z,z}<1\}$$
the complex $n$-ball in the space $\calZ$ with respect to the positive
definite form $-\spitz{\cdot,\cdot}$. There is a natural bijection
$$\calB_\calZ\Isom\calB,\quad z\loma [e+z].$$
We carry over the action of $\U(V)$ to $\calB_\calZ$ and denote it by $g\spitz z$,
$$g\spitz z:=(a+b(z))^{-1}(c+d(z)).$$
Let $g\in\GL(V)$ be an element with the property $g(e)=e$. Then $g$ acts on $V/\cz e$.
We denote by $P\subset\GL(V)$ the subgroup
$$P:=\{p\in\GL(V);\quad p(e)=e,\ p\ \hbox{acts as identity on}\ V/\cz e\}.$$
The corresponding matrices then are of the form
$$p=\pmatrix{1&b\cr0&\id_\calZ},\quad b\in \calZ^*.$$
The group $P$ is a closed complex Lie subgroup. The quotient $\GL(V)/P$ carries a natural structure as
complex manifold.
For $g\in \GL(V)$, the element $g(e)$ depends only on the coset $gP$.
Hence, the subset
$$\calB^*=\{gP\in \GL(V)/P;\quad g(e)\in\tilde\calB\}$$
is a well-defined subset of $\GL(V)/P$. It is open and hence a complex manifold too.
There are natural (holomorphic) maps
$$
\calB^*\lo\tilde\calB\lo\calB,\quad, gP\loma g(e)\loma [g(e)].$$
We consider the group
$$K_\cz=\GL(\cz e)\times\GL(\calZ)\cong\cz^*\times\GL(n,\cz)$$
as a subgroup of $\GL(V)$ in the obvious way.
The corresponding matrices are of the form
$$k=\pmatrix{k_1&0\cr 0&k_2}.$$
Usually the element $k_1$ will be identified with the corresponding complex number.
The group $K_\cz$ is the complexification of the maximal compact subgroup
$$K:=\U(\cz e)\times\U(\calZ)$$
of $\U(V)$.
\smallskip
The elements of $K_\cz$ fix the point $[e]\in \pz(V)$. Hence we have natural map
$K_\cz\to\calB^*$.
\proclaim
{Lemma}
{The natural map $K_\cz\to\calB^*$ gives a bijection between $K_\cz$ and the fibre of
the natural projection $\calB^*\to\calB$ over $[e]$.}
FibPro%
\finishproclaim
{\it Proof.\/} The elements which stabilize $[e]$ are of the form
$$g=\pmatrix{a&b\cr 0&d}.$$
They can be written in a unique way in the form $g=kp$, $k\in K_\cz$, $p\in P$.\qed
\smallskip
The group $K_\cz$ normalizes $P$ and hence acts on $G/P$ by multiplication from the right,
$$(gP,k)\loma gkP,\quad g\in\GL(V),\quad k\in K_\cz.$$
Hence $\calB^*\to \calB$ is a principal fibre bundle with structural group $K_\cz$.
\smallskip
As we mentioned already, the unitary group $\U(V)$ acts on $\tilde\calB$.
Hence it acts also on $\calB^*$ by multiplication from the left.
\smallskip
We can now  define vector valued automorphic forms.
Since $\calB^*$ plays the role of an extension of the ball $\calB$, we use from now
on letters as $z$ to denote the elements of $\calB^*$. The action of $\U(V)$ is denoted
by $\gamma z$ and that of $K_\cz$ by $zk$.
\proclaim
{Definition}
{Let $\Gamma\subset\U(V)$ a subgroup,
$\chi$ a character of $\Gamma$ and $\varrho:K_\cz\to\GL(\calU)$ a rational representation
of $K_\cz$ on some finite dimensional complex vector space. An automorphic form
for $(\Gamma,\chi,\varrho)$ is a holomorphic function
$$g:\calB^*\lo\calU$$
with the transformation property
$$f(\gamma z k)=\chi(\gamma)\varrho(k)^{-1}f(z).$$
In the case $n=1$ the usual regularity condition at the cusps has to be added.}
DefAu%
\finishproclaim
We denote the space of theses forms by $[\Gamma,\chi,\varrho]$. For trivial $\chi$ we simply write
$[\Gamma,\varrho]$.
It may happen that elements of the form $\zeta\id_V$, $\betr\zeta=1$,  are contained in $\Gamma$.
The corresponding transformations of $\calB^*$ come also from $K_\cz$. Hence
$\chi$ and $\varrho$ have to satisfy a compatibility condition if non-zero
automorphic forms exist.
\smallskip
We explain briefly the relation to the notion of
(scalar valued) automorphic form as it has been
used by Borcherds. An automorphic form in his sense is a holomorphic function
$f:\tilde\calB\to\cz$ with the transformation property
$f(\gamma z)=\chi(\gamma)f(z)$ and $f(tz)=t^{-r}f(z)$. The composition of $f$
with the projection $\calB^*\to\tilde\calB$ then gives an automorphic form in the
sense of Definition \DefAu\ with respect to the representation
$\varrho(k_1,k_2)=k_1^r$.
\smallskip
In older contexts,  automorphic forms are functions on $\calB_\calZ$ transforming with respect
to an automorphy factor. We want to describe the link between the two approaches.
For this we construct a section $\calB_\calZ\to\calB^*$.
First we construct a section $\calB\to\tilde\calB$. Each element of $V$ can be written in the
form $v=Ce+z$ where $C$ is a complex number and $z\in \calZ$.
From $\spitz{v,v}>0$ follows $C\ne 0$. Hence each element of $\calB$ has a unique representant
in $\tilde\calB$ with $C=1$. This gives a section $\calB\to\tilde\calB$.
Let now $v=Ce+w\in\tilde\calB$. We associate to $v$ a linear transformation
$g_v\in\GL(V)$, namely
$$g_v(xe+y)=Cxe+wx+y\quad(x\in\cz,\ y\in \calZ),$$
or, in matrix notation
$$g_v=\pmatrix{C&0\cr w&\id_\calZ}\qquad(v=Ce+w).$$
We have $g_v(e)=v$. Hence $g_vP$ is contained in $\calB^*$.
This gives us the desired section $\tilde\calB\to \calB^*$. Combining it with
$\calB\to\tilde\calB$ we get a section
$\calB\lo\calB^*$. Moreover  using a the isomorphism $\calB_\calZ\cong\calB$, we get
the map
$$\sigma:\calB_\calZ\lo\calB^*,\quad z\loma\pmatrix{1&0\cr z&\id_\calZ}P.$$
\proclaim
{Lemma}
{There is a ``canonical factor of automorphy''
$$J_{\hbox{\sevenrm can}}:\U(V)\times \calB_\calZ\lo K_\cz$$
with the property
$$\sigma(\gamma\spitz z)J_{\hbox{\sevenrm can}}(\gamma,z)=
\gamma\sigma(z),\qquad \gamma=\pmatrix{a&b\cr c&d}.$$
It can be defined by the formula
$$J_{\hbox{\sevenrm can}}\left(\pmatrix{a&b\cr c&d},z\right)=
\pmatrix{a+b(z)&0\cr0&d-(a+b(z))^{-1}(c+d(z))\otimes b}.$$
}
sigK%
\finishproclaim
{\it Proof.\/} We have
$$\sigma(\gamma\spitz z)=\pmatrix{1&0\cr (a+b(z))^{-1}(c+d(z))&\id}P,\qquad
\gamma\sigma(z)=\pmatrix{a+b(z)&b\cr c+d(z)&d}P.$$
The equation
$$\eqalign{
&\pmatrix{1&0\cr (a+b(z))^{-1}(c+d(z))&\id}
\pmatrix{a+b(z)&b\cr0&d-(a+b(z))^{-1}(c+d(z))\otimes b}\cr
&\qquad=\pmatrix{a+b(z)&b\cr c+d(z)&d}\cr}$$
gives the second statement of Lemma \sigK. It also implies that $J$ is an
automorphy factor.\qed
\smallskip
We call $J_{\hbox{\sevenrm can}}$ the {\it canonical automorphy factor.}
For any representation $\varrho$
of $K_\cz$ we then can define the automorphy factor
$$J_\varrho(\gamma,z)=\varrho(J_{\hbox{\sevenrm can}}(\gamma,z)).$$
If one takes for $\varrho$
the tautological representation $\id_{K_\cz}$, one obtains back the canonical
automorphy factor.
\proclaim
{Lemma}
{Let $f:\calB^*\to\calZ$ be an automorphic form with respect to $(\Gamma,\chi,\varrho)$.
Then $F(z)=f(\sigma z)$ has the transformation property
$$F(\gamma \spitz z)=\chi(\gamma)J_\varrho(\gamma,z)F(z)$$
and every holomorphic $F$ with this transformation property comes from an $f$.}
AutHa%
\finishproclaim
{\it Proof.\/} For $\gamma\in\Gamma$ we have
$$F(\gamma z)=f(\sigma\gamma\spitz z)=f(\gamma\sigma(z)J(\gamma,z)^{-1})=
v(\gamma)\varrho(J(\gamma,z))f(z).\eqno\square$$
The {\it Jacobian transformation} (derivative) $J_{\hbox{\sevenrm Jac}}(g,z)$ gives an automorphy factor of
$\U(V)$ with values in $\GL(\calZ)$.
We want to relate it to the canonical automorphy factor.
\proclaim
{Proposition}
{Consider the representation
$$\varrho: K_\cz\lo\GL(\calZ),\quad (k_1,k_2)\loma k_1^{-1}k_2.$$
(Here we consider $k_1\in\GL(\cz e)\cong\cz^*$ as complex number.)
Then
$$J_{\hbox{\sevenrm Jac}}(g,z)=J_\varrho(g,z)\quad\hbox{for}\ g\in\U(V).$$}
JacAut%
\finishproclaim
{\it Proof.\/} We will prove this not only for $g\in\U(V)$ but for all
$g\in\GL(V)$. One has to observe that both sides can be considered for arbitrary $g\in\GL(V)$
as rational functions on $\calB_\calZ$ with values in $\End(\calZ)$.
We verify the equality for generators of $\GL(V)$.
\smallni
1) $g=k=(k_1,k_2)\in K_\cz$.\hfill\break
We have
$J_{\hbox{\sevenrm can}}(k,z)=k$.
The formula
$k\spitz z=k_1^{-1}k_2 z$ shows
$$J_{\hbox{\sevenrm Jac}}(k,z)=k_1^{-1}k_2=\varrho(k)=J_\varrho(k,z).$$
\smallni
2) $g=\pmatrix{1&0\cr c&\id}$.\hfill\break
This acts as a translation $g\spitz z=z+c$ and the Jacobian is the identity.
By definition also $J_{\hbox{\sevenrm can}}(g,z)$ is the identity.
\smallni
2) $g=\pmatrix{1&b\cr 0&\id}$.\hfill\break
In this case we have
$$g\spitz z=(1+b(z))^{-1}z.$$
We have
$$J_{\hbox{\sevenrm can}}(g,z)=
\pmatrix{1+b(z)&0\cr0&\id-(1+b(z))^{-1}z\otimes b}$$
and hence
$$J_\varrho(g,z)=(1+b(z))^{-1}(\id-(1+b(z))^{-2}z\otimes b).$$
It is easy to check by means of coordinates that this is the Jacobian
of $g$.\qed
\neupara{Some examples of ball quotients}%
We consider $V=\cz^{n+1}$ and the hermitian form
$$\spitz{z,w}=\bar z_0 w_0-\bar z_1w_1-\cdots-\bar z_nw_n.$$
We denote by
$$\calE:=\gz[\zeta],\quad\zeta=e^{2\pii/3},$$
the ring of Eisenstein integers and the lattice
$$L_n=\calE^{n+1}\subset V.$$
We denote the unitary group of $L_n$ by $\Gamma_n=\U(L_n)$. We also have to consider the
congruence subgroup
$$\Gamma_n[a]=\kernel (\Gamma_n\lo\GL(n+1,\calE/a)\qquad (a\in\calE).$$
The case $a=\sqrt{-3}$ is of particular interest.
\smallskip
We are interested first in scalar valued modular forms. They belong to the
one-dimensional representation $\varrho_r(k)=k_1^r$. In the case we use the notation
$[\Gamma,\chi,r]=[\Gamma,\chi,\varrho_r]$ and we omit $\chi$ when it is trivial. The ring of
(scalar valued modular forms) is
$$A(\Gamma)=\bigoplus_{r\in\gz}[\Gamma,r].$$
The structure of this ring has been determined in the 4-dimensional case
$\Gamma_4[\sqrt{-3}]$ in [Fr] building on the paper [AF]. The corresponding modular variety
describes the variety of marked cubic surfaces. The ring $A(\Gamma_4[\sqrt{-3}])$
is rather complicated and will not be considered here. But it is possible to derive from
this 4-dimensional case several interesting cases of lower dimension. The idea is to consider
a subspace  $W\subset V$ of signature $(1,n)$, $n<4$, such that $W\cap\calE^{5}$ is a lattice
(of rank $n+1$). The embedding
$$L_{n-1}\lo L_n,\quad a\loma (a,0),$$
gives an embedding $\Gamma_{n-1}[\sqrt{-3}]\to \Gamma_n[\sqrt{-3}]$.
By restriction we obtain a ring homomorphism
$$A(\Gamma_n[\sqrt{-3}])\lo A(\Gamma_{n-1}[\sqrt{-3}]).$$
A general result states that $A(\Gamma_{n-1}[\sqrt{-3}])$ is the normalization of the image.
In this way one can prove the following result [FS1] (a different proof has been given in [Ko]).
\proclaim
{Theorem}
{The ring of modular forms
$A([\Gamma_3[\sqrt{-3}])$ is generated by six modular forms $T_1,\dots,T_6$ of weight 3
with the defining relations
$$T_1+\cdots+T_6=0,\quad T_1^3+\cdots+T_6^3=0.$$
The associated modular variety is a Segre cubic.
}
CalM%
\finishproclaim
We denote this Segre cubic by $\calS$ and by $\calR\subset \calS$ the ramification locus.
It can be described as follows. Let $\gamma\in\Gamma_3[\sqrt{-3}]$ be an element of finite
order which acts non trivially on $\calB_3$. By [ACT] it acts as a triflection
on $\calB_3$ and its fixed pint set is a so-called {\it short mirror}.
From [FS1] we can see that there is modular form of weight $5$ on $\Gamma_3[\sqrt{-3}]$
(but with non-trivial multiplier system) whose set of zeros is the union
of all short mirrors. The multiplicities are one.
In the notation of Definition 7.1 in [FS1] it is of the
form 
$$\chi:=B_1B_8B_{11}B_{13}B_{14}.$$
\proclaim
{Proposition}
{The ramification locus $\calS\subset\calB_3$ is the zero locus of a modular form $\chi$
of
weight $5$ with respect to $\Gamma_3[\sqrt{-3}]$ but with respect to a non-trivial 
multiplier system.}
RamLoc%
\finishproclaim
We are interested in vector valued modular forms with respect to the representation
$$\varrho_r\pmatrix{k_1&0\cr 0& k_2}=k_1^{r}k_2.$$
We denote the space of modular forms by $\calM(r)=[\Gamma_3[\sqrt{-3}],\varrho_r]$
The direct sum 
$$\calM=\bigoplus_{r\in\gz}[\Gamma_3[\sqrt{-3},\varrho_r]$$
is a module over
$$A=A(\Gamma_3[\sqrt{-3}]).$$
We want to determine its structure. 
\neupara{The tangent bundle of the Segre cubic}%
We study the following situation. Let $P(X_0,\cdots, X_n)$ be an irreducible
homogeneous polynomial and
$X\subset \pz^n(\cz)$ the associated
hypersurface and
$X_{\hbox{\sevenrm reg}}$ its regular locus. Let $D\subset\cz^{n-1}$ be an open domain
and let $t_0,\dots,t_n$ be holomorphic functions on $D$ without zeros and such that
$$D\lo X_{\hbox{\sevenrm reg}},\quad z\loma [t_0(z),\dots,t_0(z)]$$
is a holomorphic map on an open set of $X_{\hbox{\sevenrm reg}}$ . We want to describe the tangent space at a point $[b]\in X_{\hbox{\sevenrm reg}}$.

 The  projective tangent space $T_bX$ in $\pz^{n}(\cz)$ is 
defined by the equation
$$\sum_{i=0}^n (\partial_i P)(b)Y_i=0.$$
Here $\partial_i$ denotes the partial derivative by $X_i$. Since $X$ is a 
hypersurface, any solution of
$$\sum_{i=0}^n C_iY_i=0\quad (Y\in\hbox{inverse image of tangent space})$$
must be of the form
$$(C_0,\cdots,C_n)=\alpha (\partial_0 P)(b),\dots,(\partial_n P)(b))$$
with a constant $\alpha$
\smallskip
Now we write $b=t(z)$, $z\in D$. The tangent space  
$T_zD=\cz^{n-1}$ maps   into the space generated by the rows of
$$\pmatrix{t_0(z)&\dots&t_n(z)\cr
\partial_1t_0(z)&\dots&\partial_1t_n(z)\cr
\vdots&\vdots\cr 
\partial_{n-1}t_0(z)&\dots&\partial_{n-1}t_n(z)\cr}$$
We denote by $G_i$, $0\le i\le n$, the determinant of this matrix after cancellation
of the $i$-th column.
Hence we obtain
$$\det\pmatrix{Y_0&\cdots& Y_n\cr
t_0(z)&\dots&t_n(z)\cr
\partial_1t_0(z)&\dots&\partial_1t_n(z)\cr
\vdots&\vdots\cr 
\partial_{n-1}t_0(z)&\dots&\partial_{n-1}t_n(z)\cr}$$
or
$$\sum_{i=1}^n G_i(z)Y_i=0
\quad (Y\in\hbox{inverse image of tangent space}).$$
So we get
$$G_i(z)=f(z)\partial_iP(t(z))\quad\hbox{where}\quad f(z)\in \cz.$$
It is clear that $f(z)$ is a holomorphic function on $D$  and 
that it is non zero along the locus where the tangent map of $D\to P^n\cz$
is  injective.
\smallskip 
We want apply this to the Segre cubic $\calS$.
Therefore we have to consider $\calS$ as hypersurface in $\pz^4(\cz)$ 
(and not into $\pz^5(\cz)$ as in Theorem \CalM),
$$\calB_3\lo \calS\subset P^4\cz,\quad z\loma [T_1(z),\dots,T_5(z)].$$
The equation of $\calS$ with respect to this embedding is
$$S:=T_1^3+\cdots+T_5^3-(T_1+\cdots+T_5)^3.$$
We consider now the matrix
$$\pmatrix{T_1(z)&\dots&T_5(z)\cr
\partial_1T_1(z)&\dots&\partial_1T_5(z)\cr
\vdots&\vdots\cr 
\partial_{n-1}T_1(z)&\dots&\partial_{n-1}T_5(z)\cr}$$
and we  denote by $G_i$, $1\le i\le 5$, the determinant of this matrix after cancellation
of the $i$-th column. The consideration above shows the following result.
\proclaim
{Lemma}
{We have
$$G_i(z)=c\chi^2{\partial S\over \partial T_i}\qquad(c\in\cz).$$
}
GiSi%
\finishproclaim
{\it Proof.\/} We have shown above a formula 
$G_i(z)=f(z)(\partial S/\partial T_i)$ with a holomorphic function $f$ whose zero locus
is inside the ramification. It is easy to check that $f$ is a modular form. From Proposition
\RamLoc\ follows that up to a constant factor it is a power of $\chi$. The exponent
must be two as a weight consideration or the ramification index, studied in [FS1],  shows.\qed
\neupara{The structure theorem}%
We now can determine the structure of the
$A$-module  
$\calM=\bigoplus_{r\in\gz}[\Gamma_3[\sqrt{-3},\varrho_r]$. Recall
$A=A(\Gamma_3[\sqrt{-3}])$.
The elements $\{T_i,T_j\}$ can be considered as elements of $\calM(5)$. We consider the sub-module
$$\calN=\sum_{ij}A\{T_i,T_j\}.$$
It is sufficient to restrict to $1\le i,j\le 5$.
Our goal is to understand the structures of $\calM$ and $\calN$. First we determine defining relations of
$\calN$.
\proclaim
{Proposition}
{Defining relations for the module 
$$\calN=\sum_{1\le i,j\le 5}A\{T_i,T_j\}.$$
are
$$\leqalignno{&T_k\{T_i,T_j\}=
T_j\{T_i,T_k\}+T_i\{T_k,T_j\},\quad \{T_i,T_j\}+\{T_j,T_i\}=0&\hbox{\rm(1)}\cr
&\sum_{\nu=1}^5(\partial_\nu S)\{T_\nu,T_\mu\}=0&\hbox{\rm(2)}\cr}
$$
}
DeRel%
\finishproclaim
We recall that
$$S:=T_1^3+\cdots+T_5^3-(T_1+\cdots+T_5)^3$$
is the equation of the Segre cubic (considered as hypersurface in $\pz^4(\cz$) and
$\partial_\nu S$ denotes its derivative by $T_\nu$.
\smallni
{\it Proof of Proposition \DeRel.\/} As in section one we define a module
$$\calN'=\sum_{1\le i,j\le 5}A[T_i,T_j]$$
with symbols $[T_i,T_j]$ that satisfy the relations described in the proposition.
The is a natural homomorphism $\calN'\to\calN$ and we have to show that this
is an isomorphism. By Proposition \NisoN\ it is  sufficient that multiplication by
the variables $T_i$ and the $\partial_iS$ is injective. This can be done by means of
a computer.\qed
\smallskip
In Lemma \GiSi\ we proved
$$\det\pmatrix{T_1&\dots&T_4\cr
\partial_1T_1&\dots&\partial_1T_4\cr
\vdots&&\vdots\cr 
\partial_{n-1}T_1&\dots&\partial_{n-1}T_4\cr}=c\chi^2S_5.$$
We obtain
$$\det\pmatrix{T_1&T_1T_2&\dots&T_1T_4\cr
\partial_1T_1&T_1\partial_1T_2&\dots&T_1\partial_1T_4\cr
\vdots&&&\vdots\cr 
\partial_4T_1&T_1\partial_4T_2&\dots&T_1\partial_4T_4\cr}=
c\chi^2S_5T_1^3
$$
If we multiply the first column by $T_2$ and subtract it to the second one and so on,
we obtain the following Lemma.
\proclaim
{Lemma}
{We have
$$\det(\{T_1,T_2\},\{T_1,T_3\},\{T_1,T_4\})=c\chi^2S_5T_1^2.$$
}
KlaDet%
\finishproclaim
Since the determinant is different from $0$ every element of $\calM$ can be written in the
form
$$g_1\{T_1,T_2\}+g_2\{T_1,T_3\}+g_3\{T_1,T_4\}$$
with meromorphic functions. It is easy to check that these are meromorphic modular forms. 
In particular, they have trivial multipliers.
From Lemma \KlaDet\ we get that
the product of $h_i=g_i\chi^2S_5T_1^2$ is holomorphic. The multipliers of $\chi$ are non-trivial
on the triflections. They are third roots of unity. Hence $h_i/\chi$ is holomorphic and,
applying the same argument, $h_/\chi^2$ is holomorphic. We have shown that
$$\calM\subset {1\over T_1^2S_5}\sum_{1\le i, j\le 5}\calN.$$
During the proof we selected 1 and 5 from $\{1,\dots,5\}$. Since we could have chosen other indices
we obtain the following proposition.
\proclaim
{Proposition}
{We have
$$\calM=\bigcap_{1\le i<j\le 5}{1\over T_i^2S_j}\calN.$$
}
ProInt%
\finishproclaim
{\it Proof.} Since the elements on the right hand side are holomorphic, they must
belong to $\calM$.\qed
\smallskip
We know generators and defining relations
of $\calN$, thus  the following lemma can be proved by means of SINGULAR.
\proclaim
{Lemma}
{For arbitray $1\le i< j\le 6$ one has
$$\calN={1\over S_i}\calN\cap{1\over S_j}\calN$$}
ProIntS%
\finishproclaim
Together with Lemma \ProInt\ we obtain the following result.
\proclaim
{Proposition}
{We have
$$\calM=\bigcap_{1\le i\le 5}{1\over T_i^2}\calN.$$
}
ProIntT%
\finishproclaim
The modules $\calM$ and $\calN$ are different as the following example shows.
\proclaim
{Lemma}
{The element 
$$\eqalign{&{1\over T_1T_2}\bigl(\cr
&(2T_1T_3+T_3^2+2T_1T_4+2T_3T_4+T_4^2+2T_1T_5+2T_3T_5+2T_4T_5+T_5^2)\{T_1,T_2\}+\cr
&(-T_2T_3-2T_2T_4-2T_2T_5)\{T_1,T_3\}+\cr
&(-T_2T_4-2T_2T_5)\{T_1,T_4\}+\cr
&(-T_2T_5)\{T_1,T_5\}\bigr)\cr}
$$
is contained in $\calM$ but not in $\calN$.}
MnotN%
\finishproclaim
{\it Proof.} One has to show that that the produce of this element by $T_1T_2$
is contained in $T_1\calN\cap T_2\calN$ but not in $T_1T_2\calN$. Since we know 
the structure of $\calN$ this can be verified with the help of a computer.\qed
\smallskip
Our main result states:
\proclaim
{Theorem}
{The module $\calM$ is generated by the $\{T_i,T_j\}$ and by the special element described in
Lemma \MnotN. 
\smallni
{\bf Hilbert Functions.} (Recall that the $T_i$ have degree 3 and the $\{T_i,T_j\}$ are counted
with degree $5$.)
\smallni
The Hilbert function of $\calM$ is
$$\eqalign{&
{-t^{17} + 2t^{14} - 5t^{11} - 10t^8 - 10t^5\over (t^3-1)^3}=\cr&\quad
11t^5 + 41t^8 + 95t^{11} + 173t^{14} + 275t^{17} + 401t^{20}+\dots.
\cr}$$
The Hilbert function of $\calN$ is
$$\eqalign{&
{-5t^{11} - 8t^8 - 11t^5\over (t^3-1)^3}=\cr&\quad
10t^5 + 40t^8 + 95t^{11} + 173t^{14} + 275t^{17} + 401t^{20}+\dots.
\cr}$$
In particular the modules $\calM$ and $\calN$ differ only in the two lowest
degrees $5$ and $8$.}
MT%
\finishproclaim
\vskip1cm\noindent
{\paragratit References}%
\medskip
\item{[AF]} Allcock, D.\ Freitag, E.:
{\it Cubic Surfaces and Borcherds Products,\/}
Commentarii Math. Helv. Vol. {\bf 77}, Issue 2, 270--296 (2002)
 \medskip
\item{[CG]} Cl\'ery, F., van der Geer, G.:
{\it Generators for modules of vector-valued Picard modular forms,} 
arXiv: 1203.0131v3  (2012)
\medskip
\item{[F]} Freitag, E.: {\it A graded algebra related to cubic surfaces,\/}
Kyushu Journal of Math. Vol. {\bf 56}, No. 2, 299--312 (2002)
\medskip
\item{[FS1]} Freitag, E.,  Salvati Manni, R.:
{\it A three dimensional ball quotient,}
Math.\ Z. Vol.~{276}, Issue 1-2, pp 345-370 ( 2014) arXiv:1201.0131
\medni
\item{[FS2]} Freitag, E.,  Salvati Manni, R.:
{\it Basic vector valued Siegel modular forms of genus two,}
\medni
\item{[Ko]} Kondo, S.: {\it The Segre cubic and Borcherds products,\/}
in "Arithmetic and Geometry of K3 surfaces and Calabi-Yau threefolds", Fields Institute Communications 67, 549--565, Springer 2013
\medni
\bye